\documentclass[conference]{IEEEtran}
\IEEEoverridecommandlockouts
\usepackage{cite}
\usepackage{amsmath,amssymb,amsfonts}
\usepackage{algorithmic}
\usepackage{graphicx}
\usepackage{textcomp}
\usepackage{xcolor}
\def\BibTeX{{\rm B\kern-.05em{\sc i\kern-.025em b}\kern-.08em
    T\kern-.1667em\lower.7ex\hbox{E}\kern-.125emX}}

\begin{document}
\bstctlcite{IEEEexample:BSTcontrol}
\title{Demand Variation Impact on Tightness of Convex Relaxation Approaches for the ACOPF Problem\\
\thanks{978-1-7281-8192-9/21/\$31.00 ©2021 IEEE}
}

\author{\IEEEauthorblockN{ Arash Farokhi-Soofi}
\IEEEauthorblockA{\textit{Department of Electrical and Computer Engineering} \\
\textit{University of California San Diego }\\
San Diego, USA \\
afarokhi@ucsd.edu}
\and
\IEEEauthorblockN{ Saeed D. Manshadi}
\IEEEauthorblockA{\textit{Department of Electrical and Computer Engineering} \\
\textit{San Diego State University }\\
San Diego, USA \\
smanshadi@sdsu.edu}
}

\maketitle

\begin{abstract}
 This paper investigates the impact of the changes in the demand of power systems on the quality of the solution procured by the convex relaxation methods for the AC optimal power flow (ACOPF) problem. This investigation needs various measures to evaluate the tightness of the solution procured by the convex relaxation approaches. Therefore, three tightness measures are leveraged to illustrate the performance of convex relaxation methods under different demand scenarios. The main issue of convex relaxation methods is recovering an optimal solution which is not necessarily feasible for the original non-convex problem in networks with cycles. Thus, a cycle measure is introduced to evaluate the performance of relaxation schemes. The presented case study investigates the merit of using various tightness measures to evaluate the performance of various relaxation methods under different circumstances.
\end{abstract}


\section*{Nomenclature}

\subsection*{Parameters}
\noindent \begin{tabular}{ l p{6.55cm} }
\color{black}$C_g$ &Generation cost function of generator $g$ \\


$G_{(.)}$  & Elements of the conductance matrix\\

$B_{(.)}$  & Elements of the susceptance matrix\\

$V_i^{min},V_i^{max}$  & Min/max voltage magnitude at bus $i$\\


\color{black}$p_g^{min},p_g^{max}$&Real power generation limits of unit $g$ \\

\color{black}$q_g^{min},q_g^{max}$  &Reactive power generation limits of  unit $g$ \\

$S_{ij}^{max}$  & Maximum apparent power of the branch connecting buses $i$ and $j$\\


$p_i^d$ & Real power demand at bus $i$ \\
$q_i^d$ & Reactive power demand at bus $i$ \\





\end{tabular}

\subsection*{Variables}
\noindent \begin{tabular}{ l p{6.55cm} } 
$V_{i}$ & Voltage phasor of bus $i$ \\

$e_i$ & Real part of voltage phasor of bus $i$ \\

$f_i$ & Imaginary part of voltage phasor of bus $i$ \\

$p_g$ & Real power generation output of unit $g$ \\

$q_g$ & Reactive power generation output of unit $g$ \\


$c,s$  & Lifting operator terms for SOCP relaxation\\


$S_{ij}$ &  Apparent power flow from bus $i$ to bus $j$\\ 

$p_{ij}$ &  Real power flow from bus $i$ to bus $j$\\

$q_{ij}$ &  Reactive power flow from bus $i$ to bus $j$\\


$\gamma$ & Lifting operator term for the moment relaxation matrix\\

\end{tabular}
\subsection*{Sets}
\noindent \begin{tabular}{ l p{6.55cm} }
$\mathcal{N}$ & Set of all buses \\

$\mathcal{L}$ & Set of distribution lines  \\

$\mathcal{G}$ & Set of generation units  \\

$\mathcal{G}_i$ & Set of all generation units connected to bus $i$ \\

$\mathcal{\delta}_i$ & Set of all buses connected to bus $i$ \\
\end{tabular}
\section{Introduction}
The optimal power flow is the basic problem of many power system planning problems and its AC form is presented in \cite{dommel1968optimal}. The ACOPF is a non-linear non-convex optimization problem and it is hard to solve. Thus, the traditional iterative methods may fail to find an optimal solution for this problem. To solve this issue a set of convex relaxation approaches is introduced. Semi-definite Programming (SDP) \cite{lavaei2011zero} and Second-order Cone Programming (SOCP) \cite{jabr2007conic} relaxations are among the most popular relaxation methods. 
The major issue is that the solution procured by the relaxed problem may be infeasible for the original ACOPF problem \cite{kocuk2015inexactness}. The SOCP relaxation method is exact for radial networks. Therefore the solution procured by the relaxed problem is the same as the solution procured by the original non-convex ACOPF problem. However, the SOCP relaxation method is not exact for mesh networks due to the lack of constraints enforcing the voltage angles across the cycles \cite{kocuk2016strong}. A method to enforce the voltage angles across the cycles leveraging the second-order moment relaxation matrix of maximal cliques is presented in \cite{soofi2020socp}. Another way to improve the exactness of the SOCP relaxation method is bound tightening. Several bound tightening methods and McCormick envelops are presented in \cite{shchetinin2019efficient} and \cite{bynum2018strengthened}. To improve the exactness of the SDP relaxation method, moment relaxation based approaches are introduced in \cite{molzahn2014moment}. Employing moment-based approaches guarantees that when the order of the moment relaxation goes to infinity, the relaxation is exact. 

One way to evaluate the exactness of the solution procured by the relaxation method is through tightness measures. In the literature, two different tightness measures were introduced. The gap of the objective value with the objective value procured by solving the non-convex ACOPF problem is presented in \cite{coffrin2015qc} as a tightness measure. Another tightness measure which leverages the eigenvalues of the SDP matrix is presented in \cite{manshadi2019distributed} and \cite{manshadi2019convex}. Another way to evaluate the tightness of the SOCP relaxation method is to sum the voltage angles over a cycle that is utilized in this paper.   
The reports on the tightness of different case studies might not be accurate as they are usually presented for a selected demand profile. This paper aims to find out \textit{what is the impact of changes in the demand of the system on the tightness of convex relaxation methods?} 

\section{Problem Formulation}

The problem formulation of ACOPF in rectangular form is presented in (\ref{opf}). The objective function which is presented in (\ref{opfa}), aims to minimize the total cost of generation. The real and reactive power sending from bus $i$ to bus $j$ are presented in (\ref{opfb}) and (\ref{opfc}). The nodal real and reactive power balance of bus $i$ are given in (\ref{opfd}) and (\ref{opfe}), respectively. The upper and lower limits of voltages of buses are given in (\ref{opff}). The real and reactive generation limits are presented in (\ref{opfg}) and (\ref{opfh}), respectively. The line limit constraint is presented in (\ref{opfi}). The ACOPF problem presented in (\ref{opf}) is reformulated into a convex relaxation form as discussed in the next section.   
\begin{subequations} \label{opf}
\begin{alignat}{3}
&min\color{black}\sum_{g\in  \mathcal{G} }^{ }C_{g}(p_{{g}})\label{opfa}\\
&\text{s.t.}\nonumber\\
&\color{black}p_{ij}= G_{ij}(e_ie_j+f_if_j)-B_{ij}(e_if_j-e_jf_i)
 \hspace{0.6cm}\forall{(i,j)} \in \mathcal{L} \label{opfb}\\
&\color{black}q_{ij}=-B_{ij}(e_ie_j+f_if_j)-G_{ij}(e_if_j-e_jf_i)
\hspace{0.35cm}\forall{(i,j)} \in \mathcal{L}
\label{opfc}\\ 
&\sum_{g \in \mathcal{G}_i } p_g-p^d_i=\color{black}G_{ii}(e_i^2+f_i^2)+\sum_{j\in \mathcal{\delta}_i}p_{ij}\color{black} \hspace{1.10cm} \forall{i \in \mathcal{N}} \label{opfd}\\
&\sum_{g \in \mathcal{G}_i} q_g-q^d_i=\color{black}-B_{ii}(e_i^2+f_i^2)+\color{black}\sum_{j\in \mathcal{\delta}_i}q_{ij} \hspace{0.88cm} \forall{i \in \mathcal{N}} \label{opfe}\\
&({V_i^{min}})^2\leq  e_i^2+f_i^2\leq ({V_i^{max}})^2 \hspace{2.17cm} \forall{i \in \mathcal{N}}\label{opff}\\
&p_g^{min}\leq  p_g\leq p_g^{max} \hspace{4.03cm} \forall{g \in \mathcal{G}}\label{opfg}\\
&q_g^{min}\leq  q_g\leq q_g^{max}\hspace{4.08cm} \forall{g \in \mathcal{G}}\label{opfh}\\
& 0 \leq \sqrt{p_{ij}^2+q_{ij}^2}\leq  S_{ij}^{max}\hspace{2.8cm} \forall{(i,j) \in \mathcal{L}}\label{opfi}
\end{alignat}
\end{subequations}

\section{Solution Methodology}
\subsection{Relaxation Methods}
The original ACOPF problem presented in (\ref{opf}) is a non-convex quadratic optimization problem. The source of the non-convexity is the bi-linear terms in branch flow equations presented in (\ref{opfb}) and (\ref{opfc}). To relax this non-convex optimization problem a set of lifting variables is introduced for each method.
\subsubsection{SOCP Relaxation Method}
 A set of lifting variables is presented in (\ref{cc}) to relax the original OPF problem presented in (\ref{opf}).
 \begin{subequations} \label{cc}
\begin{alignat}{2}
&c_{ii}:=e_i^2+f_i^2= V_i^*V_i \label{cc1}\\
&c_{ij}:=e_ie_j+ f_if_j=\mathbf{Re}\{V_i^*V_j \}\label{cc2} \\
&s_{ij}:=e_if_j-f_ie_j=\mathbf{Im}\{V_i^*V_j \} \label{cc3}
\end{alignat}
\end{subequations}

Reformulating the ACOPF problem presented in (\ref{opf}) leveraging the SOCP lifting variables given in (\ref{cc}), leads to SOCP problem formulation \cite{kocuk2016strong}.
\subsubsection{SDP Relaxation Method}
Once the maximal cliques of the chordal extended graph of the network are established, the first-order moment relaxation matrix associated with each maximal clique is generated using the basis given in (\ref{basis}) and its conjugate.  
\begin{equation}
 v_c=\begin{bmatrix}
 V_i & V_j & ... & V_{\left | c \right|}
\end{bmatrix}\label{basis}
 \end{equation}
A set of SDP lifting variables is introduced in (\ref{lifting-sdp}) to relax the bilinear terms in the generated first-order moment relaxation matrix given in (\ref{sdp-matrix}).
\begin{subequations}
 \begin{alignat}{3}
 & W^{\alpha}_c=v_c^*v_c=v_c^*v_c=\begin{bmatrix}
  {V_i^*V_i} & {V_i^*V_j} & ... &{V_i^*V_{\left | c \right |}}  \\ 
  {V_j^*V_i} & {V_j^*V_j} & ... &{V_j^*V_{\left | c \right |}}  \\ 
  . & . &  & .  \\ 
  . & . &  & .  \\ 
  . & . &  & .  \\ 
  {V_{\left | c \right |}^*V_i} & {V_{\left | c \right |}^*V_j} & ... &{V_{\left | c \right |}^*V_{\left | c \right |}}  \\   
\end{bmatrix}\label{sdp-matrix}\\
 &{V_i^*V_j}\xrightarrow{lifting}{\gamma_{V_i^*V_j}}\label{lifting-sdp}\color{black}
 \end{alignat}
 \end{subequations}
Leveraging the terms in the lifted SDP relaxation matrix presented in (\ref{sdp-matrix}) to reformulate the original rectangular ACOPF problem given in (\ref{opf}) and adding a positive semi-definiteness constraint for the lifted SDP relaxation matrix leads to the SDP problem formulation\cite{lavaei2011zero}. 


\subsection{Tightness Measures}
\subsubsection{Tightness Ratio Measure}
The tightness of the SDP relaxation method can be represented as the ratio of the first eigenvalue of the SDP matrix divided to the second one when the eigenvalues of the SDP matrix sorted in decreasing order. This tightness measure is presented in (\ref{TR}). When the value of this Tightness Ratio of maximal clique $c$ ($TR_c$) goes to infinity, the SDP matrix becomes rank-1 and the SDP relaxation method is tight to the original ACOPF problem i.e. the reverse of \eqref{lifting-sdp} will hold.
\begin{equation}\label{TR}
TR_c=log(\frac{\lambda_1^c}{\lambda_2^c})
\end{equation}

In (\ref{TR}), $\lambda_1^c$, $\lambda_2^c$ are the first and the second eigenvalues of the SDP matrix associated with maximal clique $c$, if the eigenvalues of SDP matrix are sorted in decreasing order.
\subsubsection{Objective Value Measure}
Another measure to evaluate the exactness of a convex optimization method is comparing their optimality gap. The optimality gap introduced as the difference between the objective value obtained from solving the nonlinear original OPF problem with IPOPT solver \cite{wachter2006implementation} as the best known feasible point for the non-convex AC-OPF problem and the objective value obtained from solving the relaxed form of OPF problem with Mosek \cite{aps2017mosek} as conic solver \cite{coffrin2015qc}. \\
 
\begin{equation}
Gap\%=
\frac{Obj_{nonlinear}-Obj_{relaxation}}{Obj_{nonlinear}}\times 100
\end{equation}\label{optimality}
The optimality gap measure presented in (\ref{optimality}) shows the difference between the objective value of the procured solution by the convex relaxation method as the lower bound and the upper bound objective value of the solution procured by solving the non-convex original OPF problem. 

\subsubsection{Cycle Measure}
Another measure to evaluate the tightness of the solution procured by the relaxation methods is cycle measure. The main challenge is the lack of convex constraints to present the voltage angles within a cycle. Therefore a way to evaluate the tightness of the solution procured by the convex relaxation method is to calculate the summation of voltage angles of buses over cycles of the graph of the network. The difference between the calculated value and zero determines the gap of the selected cycle of the network.

\section {Results}
 In this section, several case studies are leveraged to evaluate the tightness of the SOCP and SDP relaxation methods with various tightness measures. In the first subsection the cycle measure leveraged to compare the performance of the SOCP and SDP relaxation methods with different demand scenarios. Note that $\lambda$ is the ratio of loads of the network.
\subsection{Cycle Measure}
\begin{figure}[h!]
{\includegraphics[width=\columnwidth]{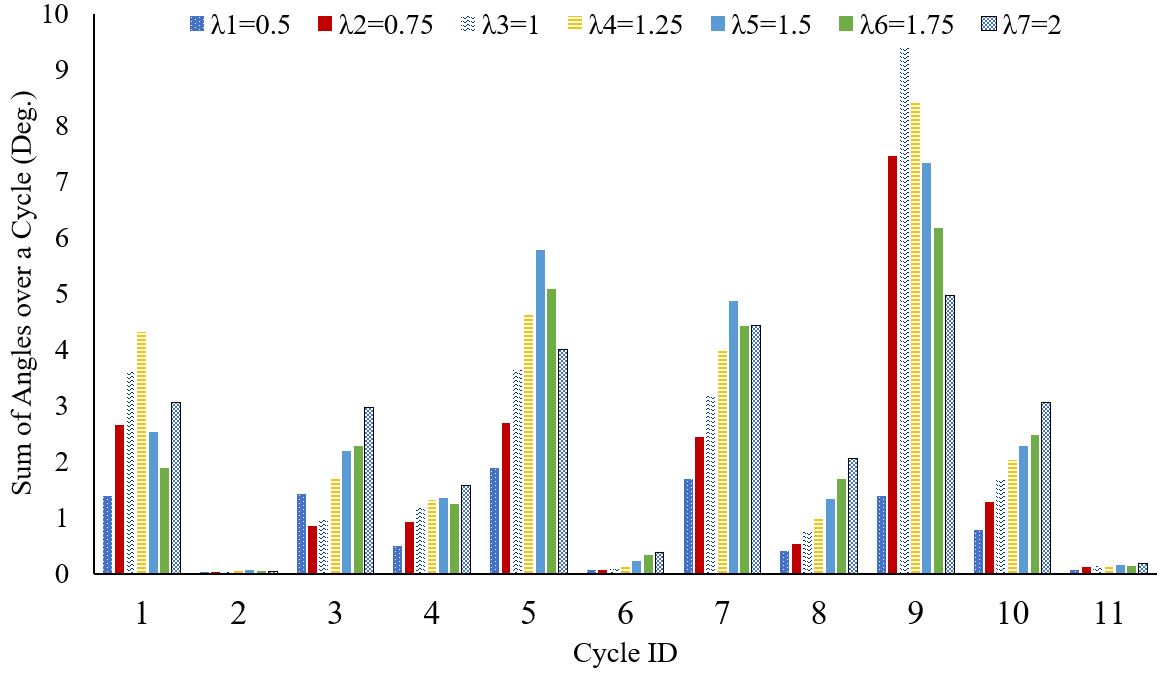}}
\caption{Profile of summation of voltage angles over a cycle procured by the SOCP method in IEEE 14-bus system}
\vspace{-0.15cm}
\label{fig:14-SOCP}
\end{figure}

\begin{figure}[h!]
{\includegraphics[width=\columnwidth]{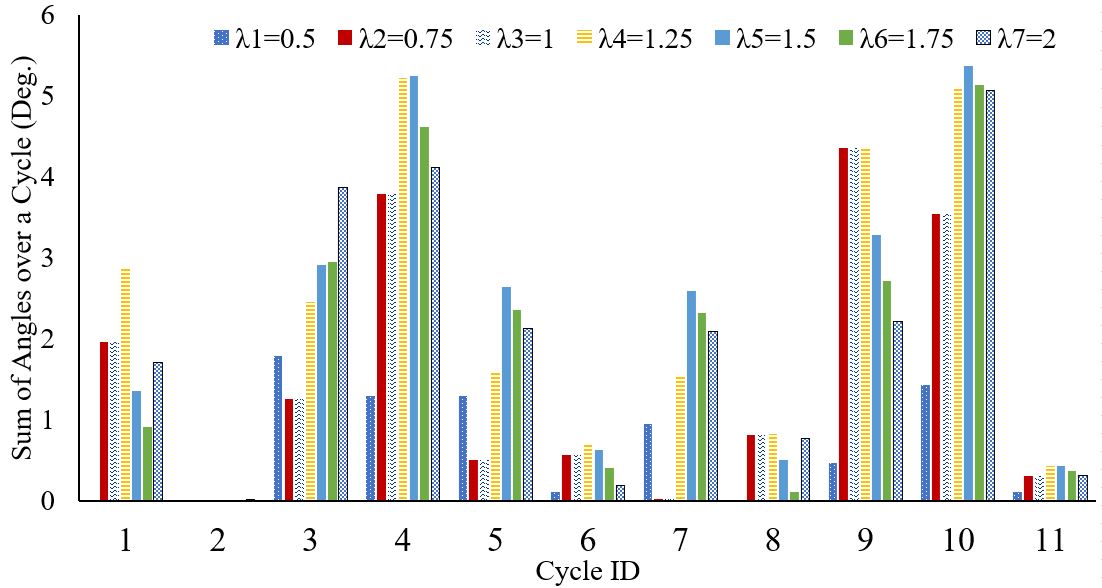}}
\caption{Profile of summation of voltage angles over a cycle procured by the SDP method in IEEE 14-bus system}
\vspace{-0.15cm}
\label{fig:14-SDP-cycle}
\end{figure}

One measure to evaluate the tightness of the relaxation methods is calculating the summation of angles of buses over every cycle of the network. In Fig. \ref{fig:14-SOCP} the cycle measure is presented for the IEEE 14-bus system when the demand of the network is changing. 
In cycle $9$, once the demand is half of the nominal value, the summation of angles obtained is decreased from $9$ degree to $2$ degree. 
In Fig. \ref{fig:14-SDP-cycle}, for the IEEE 14-bus system, the summation of voltage angles over each cycle of the network procured by the SDP relaxation method is presented.
Comparing Fig. \ref{fig:14-SOCP} and Fig. \ref{fig:14-SDP-cycle} illustrates that the average gap of cycles is less when the SDP relaxation method is employed. 

\begin{figure}[h!]
{\includegraphics[width=\columnwidth]{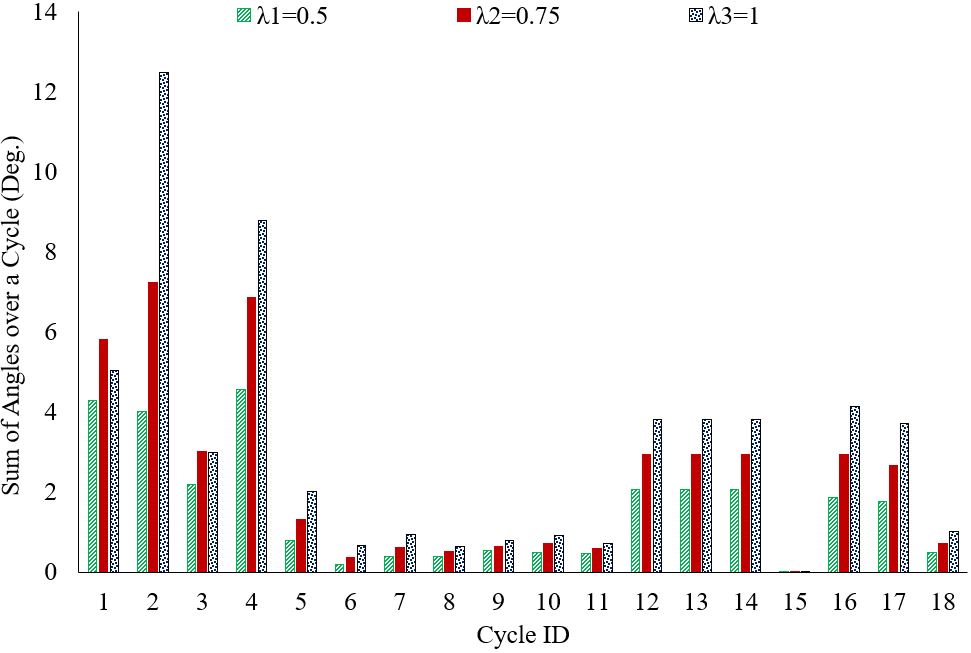}}
\caption{Profile of summation of voltage angles over a cycle procured by the SOCP method in IEEE 30-bus system}
\vspace{-0.15cm}
\label{fig:30-SOCP}
\end{figure}
\begin{figure}[h!]
{\includegraphics[width=\columnwidth]{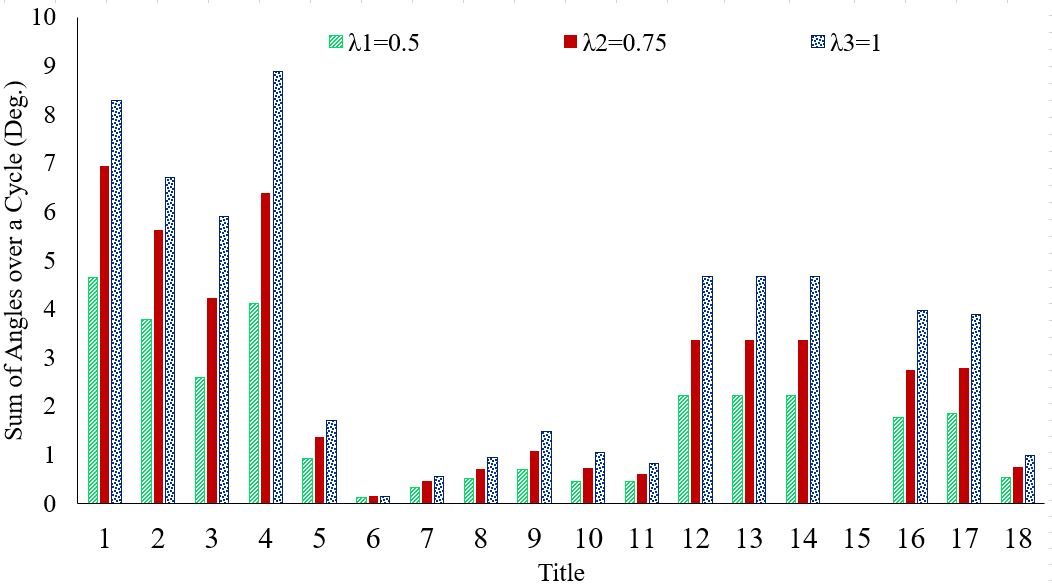}}
\caption{Profile of summation of voltage angles over a cycle procured by the SDP method in IEEE 30-bus system}
\vspace{-0.15cm}
\label{fig:30-SDP-cycle}
\end{figure}

In Fig. \ref{fig:30-SOCP}, the summation of the angles over cycles of the IEEE 30-bus system procured by the SOCP relaxation method is represented. Note that the cycles with the sum of angles more than $0.1$ degree are considered in Fig. \ref{fig:30-SOCP}. In this case, if the ratio of the demand of the network is more than $1.2$, the OPF problem is infeasible. Therefore three demand ratio is considered in this case. The largest gap of angles over a cycle occurs in cycles $2$ and $4$. In cycle $2$, when the demand ratio is half, the summation of angles over cycle decreases to $4$ degrees. In this case, when the demand is reduced to half of the nominal value, the best tightness for all cycles is obtained. In Fig. \ref{fig:30-SDP-cycle}, the summation of voltage angles over each cycle of the IEEE 30-bus system procured by the SDP relaxation method is presented. Comparing Fig. \ref{fig:30-SOCP} and Fig. \ref{fig:30-SDP-cycle} illustrates that the gap of cycles increases when the demand ratio increases. Although the cycle measure of some cycles procured by the SOCP relaxation method is less than that procured by the SDP relaxation method, the average of the sum of angles over cycles procured by the SDP relaxation method is less than the one procured by the SOCP relaxation method.

\begin{figure}[h!]
{\includegraphics[width=\columnwidth]{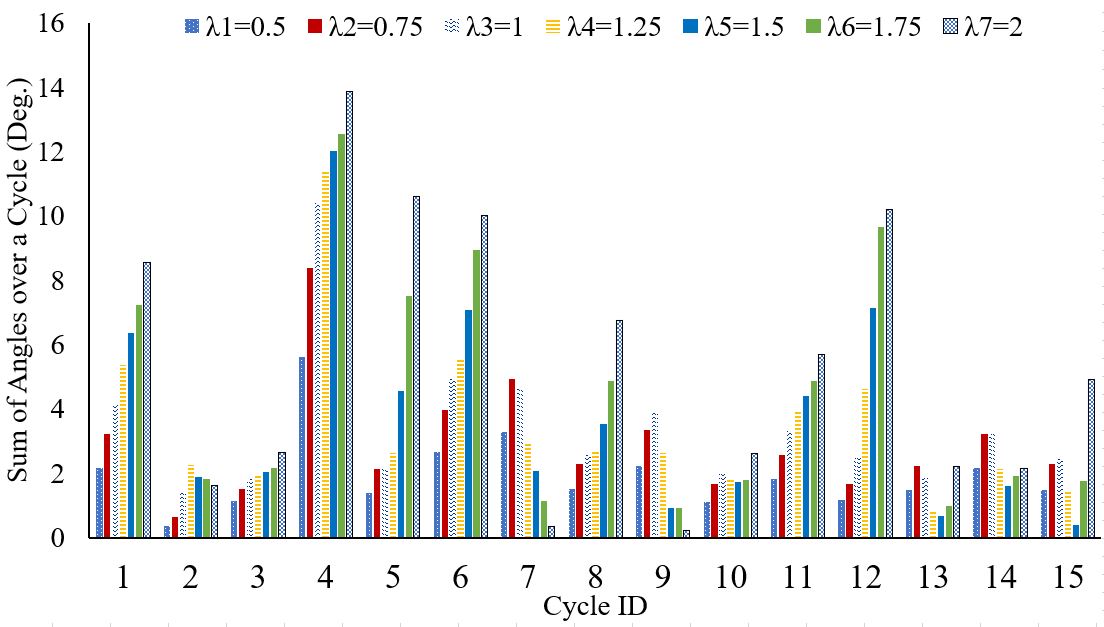}}
\caption{Profile of summation of voltage angles over a cycle procured by the SOCP method in IEEE 118-bus system}
\vspace{-0.15cm}
\label{fig:118-SOCP}
\end{figure}

\begin{figure}[h!]
{\includegraphics[width=\columnwidth]{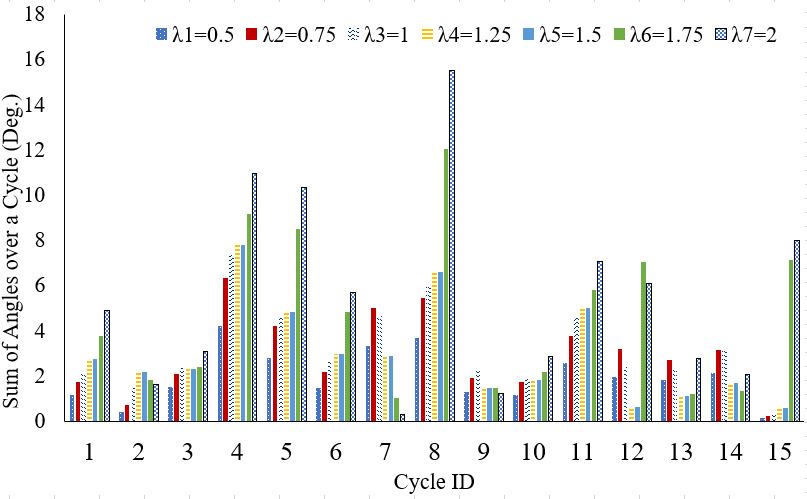}}
\caption{Profile of summation of voltage angles over a cycle procured by the SDP method in IEEE 118-bus system}
\vspace{-0.15cm}
\label{fig:118-Sdp-cycle}
\end{figure}

In Fig. \ref{fig:118-SOCP}, the summation of angles over selected cycles of the IEEE 118-bus system procured by the SOCP relaxation method is presented. The tightness of cycles when the demand ratio is $2$ is the least one at most cycles. 
The sum of angles over cycles procured by the SDP relaxation method is shown in Fig. \ref{fig:118-Sdp-cycle}. Comparing it with the one procured by the SOCP relaxation method reveals that the pattern of tightness is similar to various demand scenarios while the magnitude is higher for the SOCP method. 

\begin{figure}[h!]
{\includegraphics[width=\columnwidth]{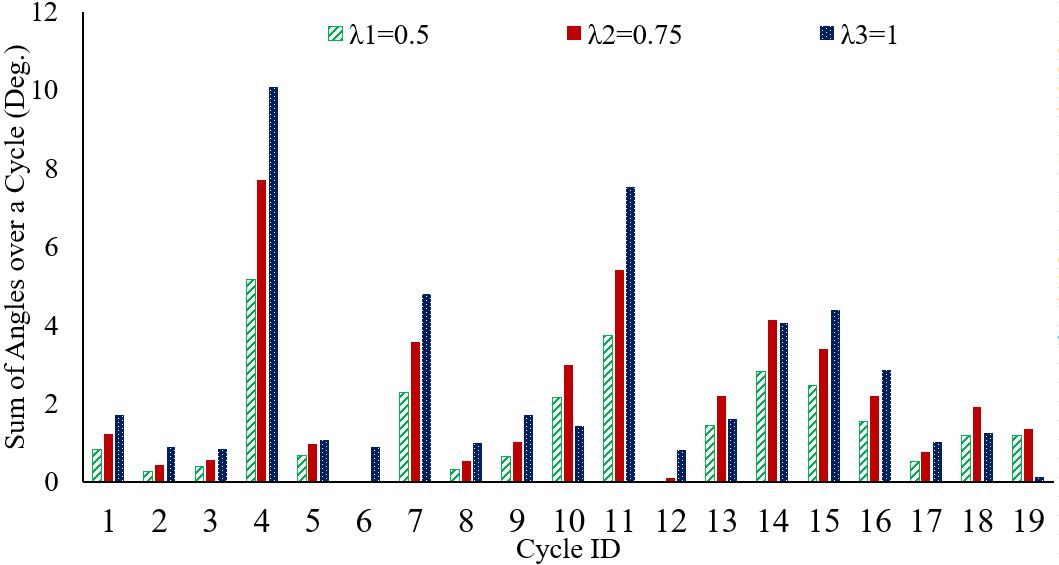}}
\caption{Profile of summation of voltage angles over a cycle procured by the SOCP method in IEEE 300-bus system}
\vspace{-0.15cm}
\label{fig:300-SOCP}
\end{figure}

\begin{figure}[h!]
{\includegraphics[width=\columnwidth]{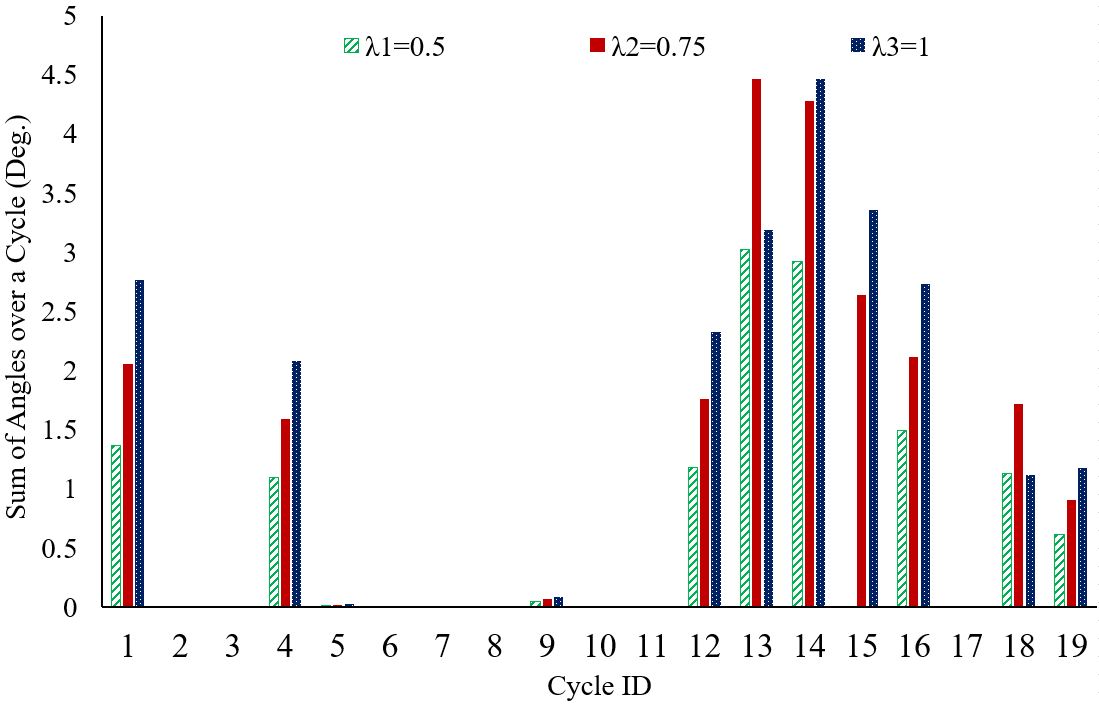}}
\caption{Profile of summation of voltage angles over a cycle procured by the SDP method in IEEE 300-bus system}
\vspace{-0.15cm}
\label{fig:300-SDP-cycle}
\end{figure}

In Figs. \ref{fig:300-SOCP} and \ref{fig:300-SDP-cycle}, the summation of angles over selected cycles of the IEEE 300-bus system procured by the SOCP and SDP relaxation methods are presented, respectively. In this case, if the ratio of the demand of the network is more than $1.1$, the OPF problem is infeasible. Note that in Fig. \ref{fig:300-SOCP}, the presented cycles are selected randomly among the cycles in which the summation of their buses voltage angles are more than $1$ degree for all demand ratios, and the same cycles are chosen to evaluate the SDP relaxation method by the cycle measure. Figs. \ref{fig:300-SOCP} and \ref{fig:300-SDP-cycle} show that the summation of angles became zero for several cycles (e.g. $2,3,6,7,8,10,11,17$) by the SDP method while compared to SOCP, it remained in a similar range for others and is increased for a few cycles (e.g. $1,12$).  

\subsection{TR Measure}
In this section, the Tightness Ratio (TR) measure evaluates the tightness of solutions procured by the SOCP and SDP relaxation methods of four test cases under different demand ratios.
\begin{figure}[h!]
{\includegraphics[width=\columnwidth]{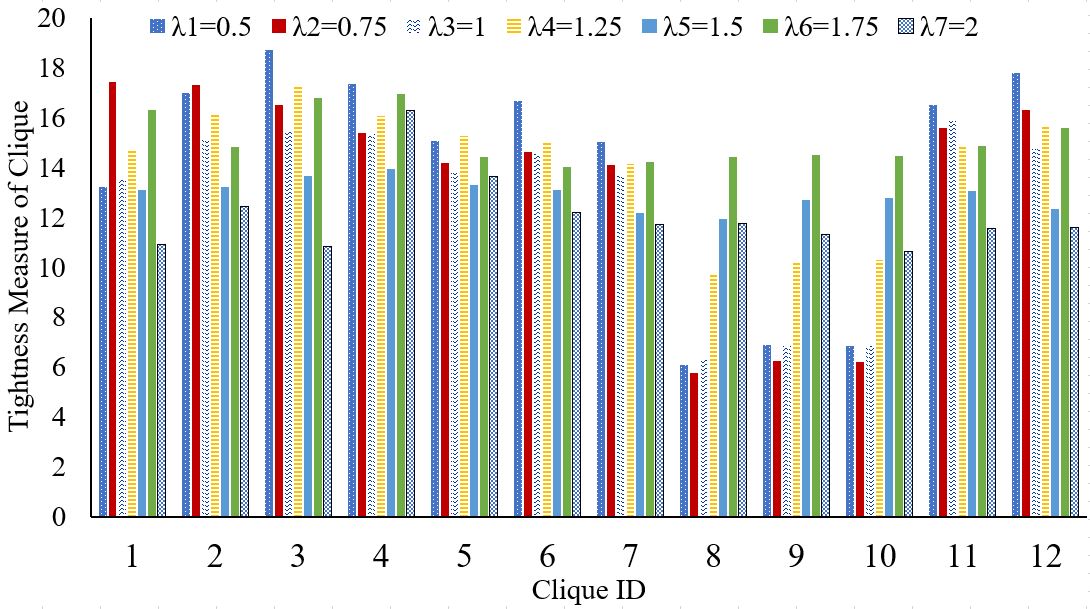}}
\caption{Profile of tightness of each clique procured by the SDP method in IEEE 14-bus system}
\vspace{-0.15cm}
\label{fig:14-SDP}
\end{figure}

\begin{figure}[h!]
{\includegraphics[width=\columnwidth]{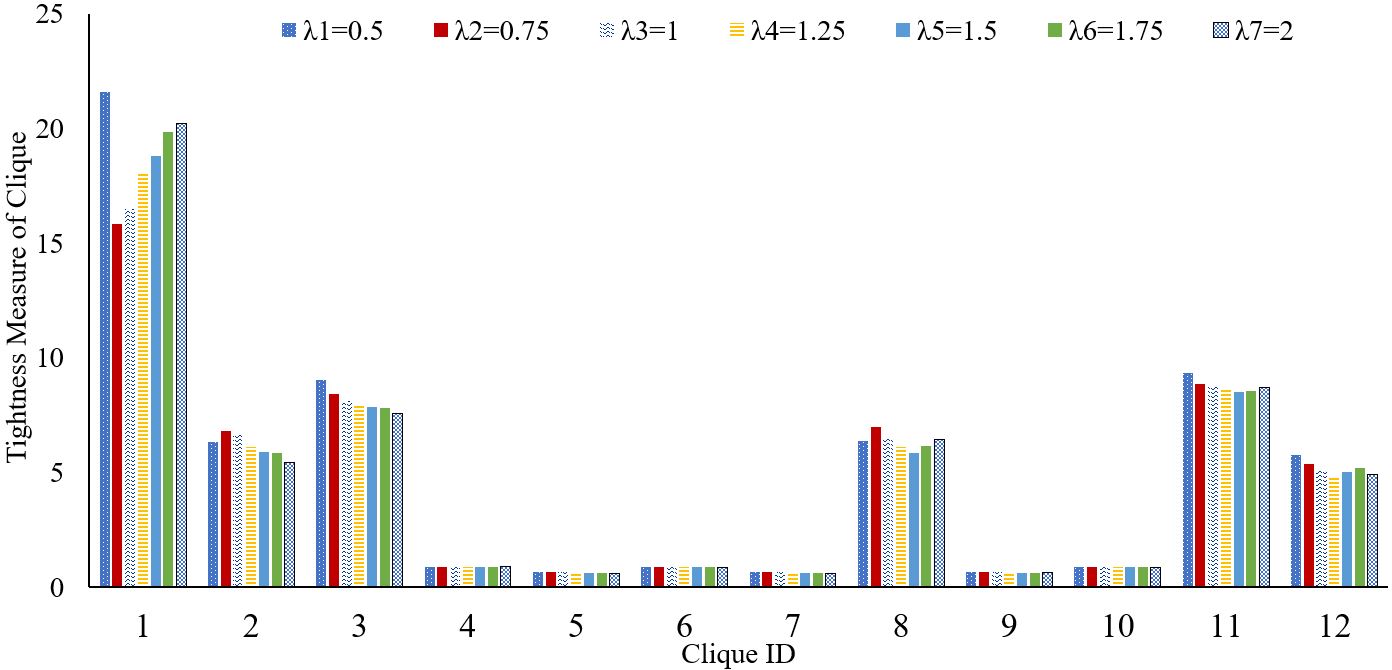}}
\caption{Profile of tightness of each clique procured by the SOCP method in IEEE 14-bus system}
\vspace{-0.15cm}
\label{fig:14-SOCP-TR}
\end{figure}
In Fig. \ref{fig:14-SDP}, the tightness ratio of each maximal clique of the network procured by the SDP relaxation method is presented. The value of TR is more than $10$ for most of the cliques under different demand ratios. This means that the first eigenvalue of the SDP matrix associated with the maximal clique is $1e^{10}$ times the second eigenvalue. Therefore the SDP relaxation method is tight. However, this tightness measure is around $6$ for cliques $8,9,10$ when the ratio of demand is $0.5,0.75,1$, respectively. Comparing the objective value gap of the SDP relaxation method in different load ratios shows that when the demand ratio is $1.25$ the optimality gap is $\%0.0001$ and the largest gap occurs when the demand ratio is $2$ and the optimality gap is $\% 0.1155$. Since the objective value gap measure considers the whole of the network and the TR measure consider each clique, this difference in finding the tightest solution of different demand ratios exists. In Fig. \ref{fig:14-SOCP-TR}, the TR measure of each maximal clique of the network procured by the SOCP relaxation method is presented. Comparing the TR of the maximal cliques procured by the SDP relaxation method with the one procured by the SOCP relaxation method illustrates that the SDP relaxation method procures a solution with less optimality gap. This conclusion is consistent with the one concluded from evaluating the SDP and SOCP relaxation methods with cycle measures in the last subsection.

\begin{figure}[h!]
{\includegraphics[width=\columnwidth]{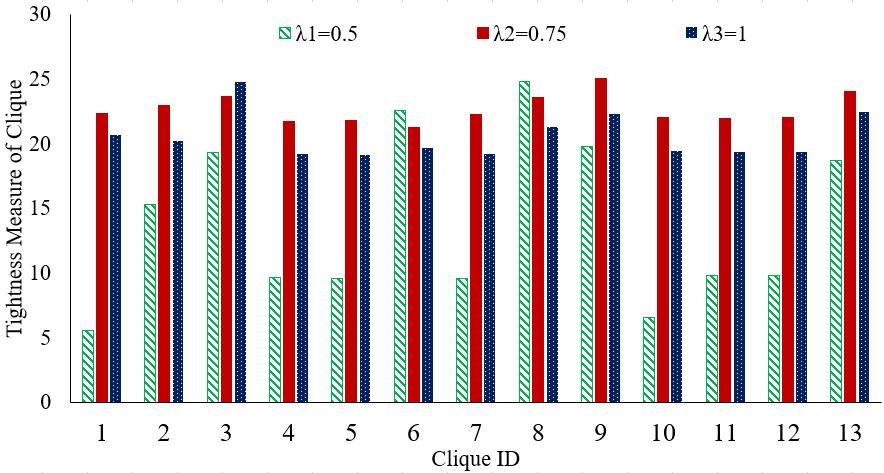}}
\caption{Profile of tightness of each clique procured by the SDP method in IEEE 30-bus system}
\vspace{-0.15cm}
\label{fig:30-SDP}
\end{figure}

\begin{figure}[h!]
{\includegraphics[width=\columnwidth]{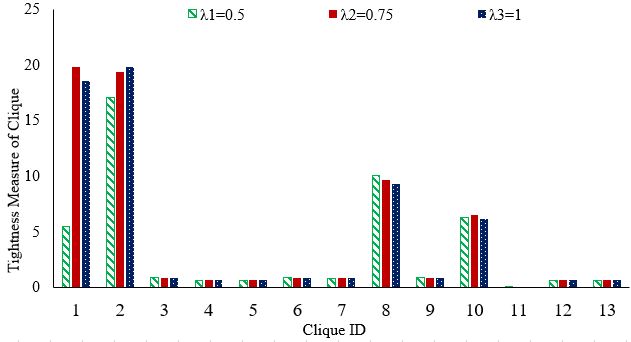}}
\caption{Profile of tightness of each clique procured by the SOCP method in IEEE 30-bus system}
\vspace{-0.15cm}
\label{fig:30-SOCP-TR}
\end{figure}

In Figs. \ref{fig:30-SDP} and \ref{fig:30-SOCP-TR}, the tightness ratio of maximal cliques of the graph associated with the IEEE 30-bus system procured by the SDP and SOCP relaxation methods are shown, respectively. Note that the cliques with tightness ratio of more than 20 for all demand ratios procured by the SDP relaxation method are not considered in Fig. \ref{fig:30-SDP} and those selected maximal cliques are chosen to evaluate the SOCP relaxation method by the TR measure. 
It is shown that with the decrease in the demand, the TR is decreasing for both SOCP and SDP relaxation for most of the cliques except clique $8$. It is also interesting that the TR is generally much smaller for the SOCP relaxation. 

\begin{figure}[h!]
{\includegraphics[width=\columnwidth]{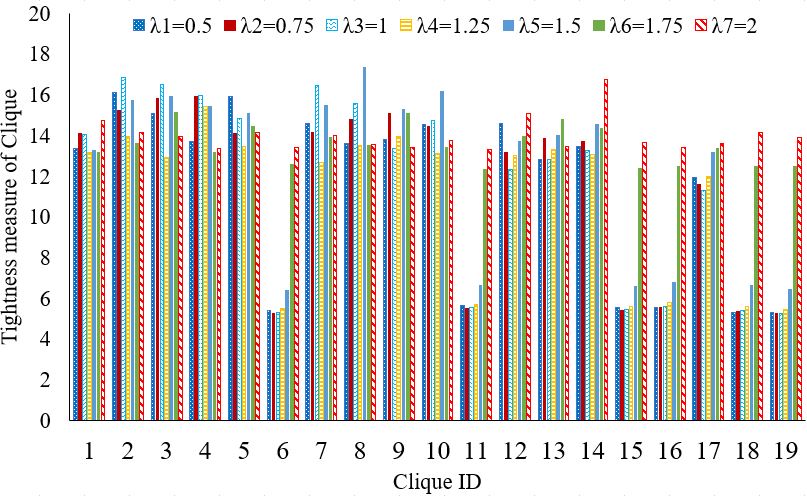}}
\caption{Profile of tightness of each clique procured by the SDP method in IEEE 118-bus system}
\vspace{-0.15cm}
\label{fig:118-SDP}
\end{figure}

\begin{figure}[h!]
{\includegraphics[width=\columnwidth]{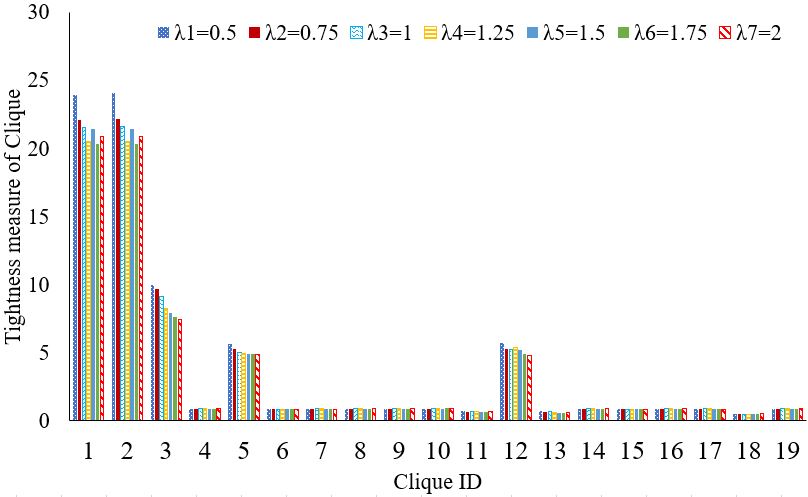}}
\caption{Profile of tightness of each clique procured by the SOCP method in IEEE 118-bus system}
\vspace{-0.15cm}
\label{fig:118-SOCP-TR}
\end{figure}
In Figs. \ref{fig:118-SDP} and \ref{fig:118-SOCP-TR}, the tightness ratio of maximal cliques of the graph associated with the IEEE 118-bus system procured by the SDP and SOCP relaxation methods are shown, respectively. An interesting observation in Fig. \ref{fig:118-SDP} is that in cliques $6,11,15,18,19$ the tightness ratio is around $6$ for demand ratios $0.5-1.5$, while the tightness ratio for these maximal cliques is increases to $12$ when the demand ratio is $1.75,2$. This means that the SDP relaxation method procured a tighter solution in higher demands for the IEEE 118-bus system. An interesting observation is the independence of the TR of the procured solution for both approaches from the utilized demand.

\begin{figure}[h!]
{\includegraphics[width=\columnwidth]{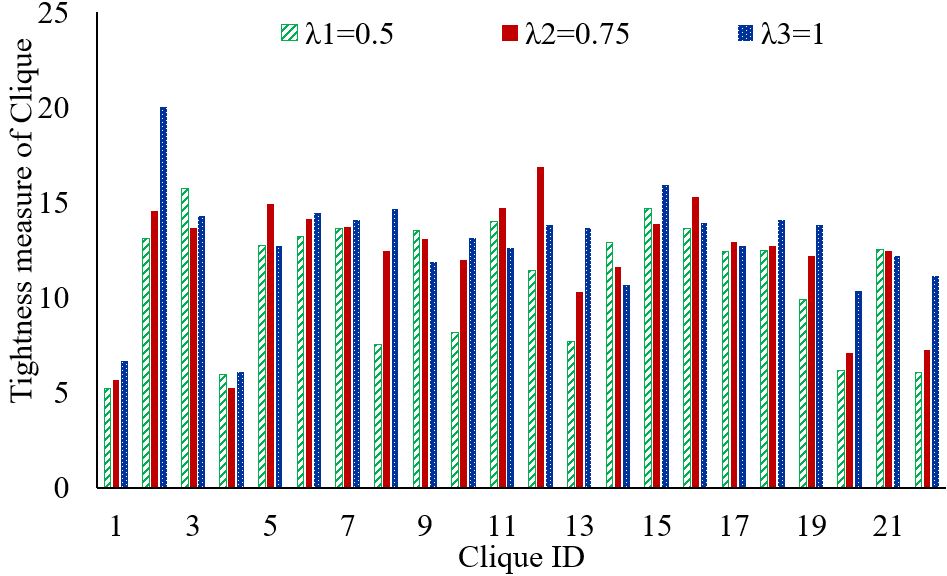}}
\caption{Profile of tightness of each clique procured by the SDP method in IEEE 300-bus system}
\vspace{-0.15cm}
\label{fig:300-SDP}
\end{figure}

\begin{figure}[h!]
{\includegraphics[width=\columnwidth]{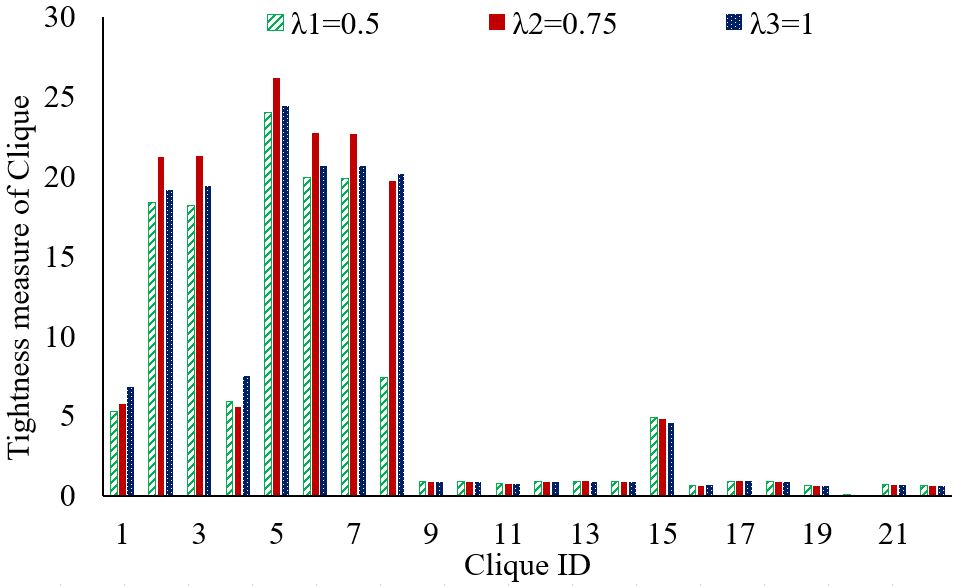}}
\caption{Profile of tightness of each clique procured by the SOCP method in IEEE 300-bus system}
\vspace{-0.15cm}
\label{fig:300-SOCP-TR}
\end{figure}
In Figs. \ref{fig:300-SDP} and \ref{fig:300-SOCP-TR}, the tightness ratio of maximal cliques of the graph associated with the IEEE 300-bus system procured by the SDP and SOCP relaxation methods are shown, respectively. Similar to the previous case studies, SDP relaxation renders a tighter solution compared to the SOCP method. 
\vspace{-0.15cm}
\section{Conclusion}
Changing the demand scenarios of the network changes the solution to the OPF problem and the relaxation gap of relaxation methods. 
Test cases elaborate that the tightness of the relaxation method is a function of the utilized demand and single demand point verification is not sufficient to rely on an approach for a given system. Another point that is illustrated in the test cases is that comparing the optimality gap calculated for each relaxation approach is not an indicator to determine if the solution that procured by each method is feasible for the original non-convex problem.
\\
\\
\\

\bibliographystyle{IEEEtran}
\bibliography{mendeley.bbl}
\end{document}